# ECONOMETRIC ANALYSIS ON EFFICIENCY OF ESTIMATOR

*M. Khoshnevisan, Lecturer, Griffith University, School of Accounting and Finance, Australia*
*F. Kaymram, Assistant Professor, Massachusetts Institute of Technology, Department of Mechanical Engineering, USA; currently at Sharif University, Iran.*
*Housila P. Singh, Rajesh Singh, Professors of Statistics, Vikram University, Department of Mathematics and Statistics, India.*
*F. Smarandache, Associate Professor, Department of Mathematics, University of New Mexico, USA.*

## ABSTRACT

This paper investigates the efficiency of an alternative to ratio estimator under the super population model with uncorrelated errors and a gamma-distributed auxiliary variable. Comparisons with usual ratio and unbiased estimators are also made.

Key words: Bias, Mean Square Error, Ratio Estimator Super Population.

## 1. INTRODUCTION

It is well known that the ratio method of estimation occupies an important place in sample surveys. When the study variate y and the auxiliary variate x is positively (high) correlated, the ratio method of estimation is quite effective in estimating the population mean of the study variate y utilizing the information on auxiliary variate x.

Consider a finite population with N units and let $x_i$ and $y_i$ denote the values for two positively correlated variates x and y respectively for the ith unit in this population, i=1,2,…,N. Assume that the population mean $\bar{X}$ of x is known. Let $\bar{x}$ and $\bar{y}$ be the sample means of x and y respectively based on a simple random sample of size n (n < N) units drawn without replacement scheme. Then the classical ratio estimator for $\bar{Y}$ is defined by

$$\bar{y}_r = \bar{y}(\bar{X}/\bar{x}) \qquad (1.1)$$

The bias and mean square error (MSE) of $\bar{y}_r$ are, up to second order moments,

$$B(\bar{y}_r) = \lambda (R S^2_x - S_{yx})/\bar{X} \qquad (1.2)$$
$$M(\bar{y}_r) = \lambda (S^2_y + R^2 S^2_x - 2R S_{yx}), \qquad (1.3)$$

where $\lambda = (N-n)/(nN)$,

$R = \bar{Y}/\bar{X}$, $S_y^2 = (N-1)^{-1}\sum_{i=1}^{N}(y_i - \bar{Y})^2$, $s_x^2 = (N-1)^{-1}\sum_{i=1}^{N}(x_i - \bar{X})^2$,

and $S_{yx} = (N-1)^{-1}\sum_{i=1}^{N}(y_i - \bar{Y})(x_i - \bar{X})$.

It is clear from (1.3) that $M(\bar{y}_r)$ will be minimum when

$$R = S_{yx}/S_x^2 = \beta, \qquad (1.4)$$

where $\beta$ is the regression coefficient of y on x. Also for $R = \beta$,
the bias of $\bar{y}_r$ in (1.2) is zero. That is, $\bar{y}_r$ is almost unbiased for $\bar{Y}$.

Let $E(\bar{y}|\bar{x}) = \alpha + \beta \bar{x}$ be the line of regression of $\bar{y}$ on $\bar{x}$, where E denotes averaging over all possible sample design simple random sampling without replacement (SRSWOR). Then $\beta = S_{yx}/S_x^2$ and $\bar{Y} = \alpha + \beta \bar{X}$ so that, in general,

$$R = (\alpha/\bar{X}) + \beta \qquad (1.5)$$

It is obvious from (1.4) and (1.5) that any transformation that brings the ratio of population means closer to $\beta$ will be helpful in reducing the mean square error (MSE) as well as the bias of the ratio estimator $\bar{y}_r$. This led Srivenkataramana and Tracy (1986) to suggest an alternative to ratio estimator $\bar{y}_r$ as

$$\bar{y}_a = \bar{z}(\bar{X}/\bar{x}) + A = \bar{y}_r - A\{(\bar{X}/\bar{x}) - 1\} \qquad (1.6)$$

which is based on the transformation

$$\bar{z} = \bar{y} - A, \qquad (1.7)$$

where $E(\bar{z}) = \bar{Z}(= \bar{Y} - A)$ and A is a suitably chosen scalar.

In this paper exact expressions of bias and MSE of $\bar{y}_a$ are worked out under a super population model and compared with the usual ratio estimator.

## 2. THE SUPER POPULATION MODEL

Following Durbin (1959) and Rao (1968) it is assumed that the finite population under consideration is itself a random sample from a super population and the relation between x and y is of the form:

$$y_i = \alpha + \beta\, x_i + u_i \; ; \quad (i = 1,2,\ldots,N)$$

where $\alpha$ and $\beta$ are unknown real constants; $u_i$'s are uncorrelated random errors with conditional (given $x_i$) expectations

$$E(u_i | x_i) = 0$$

$$E(u_i^2 | x_i) = \delta\, x_i^g$$

( i=1,2,….,N), $o \langle \delta \langle \infty$, $o \leq g \leq 2$ and $x_i$ are independently identically distributed ( i.i.d.) with a common gamma density

$$G(\theta) = e^{-x} x^{\theta-1} / \Gamma\theta,\ x \rangle o,\ 2 \langle \theta \langle \infty \ . \tag{2.1}$$

We will write $E_x$ to denote expectation operator with respect to the common distribution of $x_i$ (i=1,2,3,…,N) and $E_x E_c$, as the over all expectation operator for the model. We denote a design by p and the design expectation $E_p$, for instance, see Chaudhuri and Adhikary (1983,89) and Shah and Gupta (1987). Let 's' denote a simple random sample of N distinict labels chosen without replacement out of i=1,2,3……N. Then

$$X(=N\bar{X}) = \sum_{i \in s} x_i + \sum_{i \notin s} x_i$$

Following Rao and Webster (1966) we will utilize the distributional properties of $x_j / x_i$, $\sum_{i \in s} x_i$, $\sum_{i \notin s} x_i$, $\sum_{i \in s} x_i / \sum_{i \notin s} x_i$ in our subsequent derivations.

3. THE BIAS AND MEAN SQUARE ERROR

The estimator $\bar{y}_a$ in (1.6) can be written as

$$\bar{y}_a = \left[ (1/n)\left(\sum_{i \in s} y_i\right) \frac{\left(n\sum_{i=1}^{N} x_i\right)}{\left(N\sum_{i \in s} x_i\right)} - A\left\{ \frac{\left(n\sum_{i=1}^{N} x_i\right)}{\left(N\sum_{i \in s} x_i\right)} - 1 \right\} \right] \quad (3.1)$$

based on a simple random sample of n distinct labels chosen without replacement out of $i = 1, 2, \ldots, N$.

The bias
$$B = E_p(\bar{y}_a - \bar{Y}) \quad (3.2)$$

of $\bar{y}_a$ has model expectation $E_m(B)$ which works out as follows:

$$E_m(B(\bar{y}_a)) = E_p E_x E_c \left[ \left\{ \alpha + \beta(1/n)\left(\sum_{i \in s} x_i\right) + \bar{u} \right\} \frac{n\sum_{i=1}^{N} x_i}{n\sum_{i \in s} x_i} \right.$$

$$\left. - A\left\{ \frac{n\left(\sum_{i=1}^{N} x_i\right) - 1}{N\left(\sum_{i \in s} x_i\right)} \right\} - \right.$$

$$- E_x E_c (\alpha + \beta \bar{x} + \bar{U})$$

$$= E_p E_x E_c$$
$$\left[ \alpha\left(n\sum_{i=1}^{N} x_i / N\sum_{i \in s} x_i\right) + \beta(1/N)\left(\sum_{i=1}^{N} x_i\right) + \left(\sum_{i \in s} u_i\right)\left(\sum_{i=1}^{N} x_i / N\sum_{i \in s} x_i\right) - A\left\{\left(\left(n\sum_{i=1}^{N} x_i\right) \middle/ \left(N\sum_{i \in s} x_i\right)\right) - 1\right\} \right]$$

$$- E_x E_c (\alpha + \beta \bar{X})$$

$$= E_p E_x \left[ \alpha\left(n\sum_{i=1}^{N} x_i\right) \middle/ \left(N\sum_{i \in s} x_i\right) + \beta\bar{X} - A\left\{\left(n\sum_{i=1}^{N} x_i\right) \middle/ \left(N\sum_{i \in s} x_i\right) - 1\right\} \right] - \alpha - \beta E_x(\bar{X})$$

$$= E_x \left[ \alpha(n/N)\left(1 + \sum_{i \notin s} x_i / \sum_{i \in s} x_i\right) - A\left\{(n/N)\left(1 + \sum_{i \notin s} x_i / \sum_{i \in s} x_i\right) - 1\right\} \right] - \alpha$$

$$= \alpha(n/N)\{1 + (N-n)\theta/(n\theta - 1)\}$$
$$\qquad -A\{(n/N)(1 + (N-n)\theta/(n\theta - 1)) - 1\} - \alpha$$

$$= \alpha[(n/N - 1) + \{n(N-n)\theta / N(n\theta - 1)\}]$$

$$-A\left[-(N-n)/N+\{(N-n)n\theta/N(n\theta-1)\}\right]$$

$$= (N-n)(\alpha-A)/N(n\theta-1) \tag{3.3}$$

For SRSWOR sampling scheme, the mean square error

$$M(\bar{y}_a) = E_p(\bar{y}_a - \bar{Y})^2 \tag{3.4}$$

of $\bar{y}_a$ has the following formula for model expectations

$E_m(M(\bar{y}_a))$ :

$$E_m(M(\bar{y}_a)) = \left[E_m(M(\bar{y}_r)) + (N-n)(Nn\theta + 2N - 2n)(A^2 - 2A\alpha)/N^2(n\theta-1)(n\theta-2)\right] \tag{3.5}$$

where
$$M(\bar{y}_r) = E_p(\bar{y}_r - \bar{Y})^2 \tag{3.6}$$

is the MSE of $\bar{y}_r$ under SRSWOR scheme has the model expectation

$$E_m(M(\bar{y}_r)) = \{(N-n)/N^2\}$$

$$\left[\left\{\frac{(Nn\theta+2N-2n)\alpha^2}{(n\theta-1)(n\theta-2)}\right\} + \frac{\delta\left\{(n\theta+g-1)(n\theta+g-2)+n\theta(N\theta-n\theta+1)\right\}}{(n\theta+g-1)(n\theta+g-2)} \frac{\Gamma(\theta+g)}{\Gamma\theta}\right] \tag{3.7}$$

$[See, Rao(1968, p.439)]$

Further, we note that for SRSWOR sampling scheme, the bias

$$B(\bar{y}_r) = E_p(\bar{y}_r - \bar{Y}) \tag{3.8}$$

of usual ratio estimator has the model expectation

$$E_m(B(\bar{y}_r)) = (N-n)\alpha/(n\theta-1) \tag{3.9}$$

We note from (3.3) and (3.9) that

$$|E_m(B(\bar{y}_a))| \langle |E_m(B(\bar{y}_r))|$$

if

$$|(\alpha - A)| \langle |\alpha|$$

or if

$$(\alpha - A)^2 \langle \alpha^2$$

or if

$$o \langle A \langle 2\alpha \qquad (3.10)$$

Further we have from (3.5) that

$$E_m(M(\bar{y}_a)) - E_m(M(\bar{y}_r)) < o$$

if

$$(A^2 - 2A\alpha) < o$$

or if

$$o \langle A \langle 2\alpha \qquad (3.11)$$

which is the same as in (3.10).
Thus we state the following theorem:
Theorem 3.1 : The estimator $\bar{y}_a$ is less biased as well as more efficient than usual ratio estimator $\bar{y}_r$ if

$$o \langle A \langle 2\alpha \qquad (\alpha \neq o)$$

i . e . when A lies between $o$ and $2\alpha$.
Therefore, when intercept term $\alpha(\neq o)$ in the model (2.1) is sizable, there will be sufficient flexibility in picking A.

It is to be noted that for $\alpha = o, \bar{y}_r$ is unbiased and efficient than $\bar{y}_a$.
The minimization of (3.5) with respect to A leads to

$$A = \alpha = A_{opt} \text{ (say)} \qquad (3.12)$$

Substitution of (3.12) in (3.5) yields the minimum value of

$E_m(M(\bar{y}_a))$ as

$$\text{min. } E_m(M(\bar{y}_a)) = \frac{(N-1)}{N^2} \frac{\delta[(n\theta + g - 1)(n\theta + g - 2) + n\theta(N\theta - n\theta + 1)]}{(n\theta + g - 1)(n\theta + g - 2)} \frac{\Gamma(\theta + g)}{\Gamma \theta} \qquad (3.13)$$

which equals to $E_m(M(\bar{y}_r))$ when $\alpha = o$.
It is interesting to note that when $A = \alpha, \bar{y}_a$ is unbiased and attained its minimum average MSE in model (2.1).

In practice the value of $\alpha$ will have to be assessed, at the estimation stage, to be used as A. To assess $\alpha$, we may use scatter diagram of y versus x for data from a pilot study, or a part of the data from the actual study and judge the y-intercept of the best fitting line.

From (3.7) and (3.13) we have

$$E_m(M(\bar{y}_r)) - \min.E_m(M(\bar{y}_a)) = \{(N-n)(Nn\theta + 2N - 2n)\alpha^2\}/\{N^2(n\theta-1)(n\theta-2)\}$$
$$\rangle \; o \quad (3.14)$$

which shows that $\bar{y}_a$ is more efficient than ratio estimator when $A = \alpha$ is known exactly. For $\alpha = o$

$$\min.E_m(M(\bar{y}_a)) = E_m(M(\bar{y}_r)) \quad (3.15)$$

For SRSWOR, the variance

$$V(\bar{y}) = E_p(\bar{y} - \bar{Y})^2 \quad (3.16)$$

of usual unbiased estimator has the model expectation:

$$E_m(V(\bar{y})) = (N-n)[\beta^2\theta + \{\delta\Gamma(\theta+g)/\Gamma\theta\}]/nN \quad (3.17)$$

The expressions of $E_m(M(\bar{y}_a))$ and $E_m(V(\bar{y}))$ are not easy task to compare algebraically. Therefore in order to facilitate the comparison, denoting

$$E_1 = 100 E_m(V(\bar{y}))/E_m(M(\bar{y}_a)) \text{ and } E_2 = 100 E_m(V(\bar{y}_r))/E_m(M(\bar{y}_a)),$$

we present below in tables 1,2,3, the values of the relative efficiencies of $\bar{y}_a$ with respect to $\bar{y}$ and $\bar{y}_r$ for a few combination of the parametric values under the model (2.1). Values are given for $N = 60$, $\delta = 2.0, \theta = 8, \alpha = 0.5$, 1.0, 1.5, $\beta = 0.5, 1.0, 1.5$ and $g = 0.0, 0.5, 1.0, 1.5, 2.0$. The ranges of A, for $\bar{y}_a$ to be better than $\bar{y}_r$ for given $\alpha = 0.5, 1.0, 1.5$ are respectively (0,1), (0,2), (0,3). This clearly indicates that as the size of $\alpha$ increases the range of A for $\bar{y}_a$ to be better than $\bar{y}_r$ increases i.e. flexibility of choosing A increases.

We have made the following observations from the tables 1,2 and 3 :

(i) As g increases both $E_1$ and $E_2$ decrease. When n increases $E_1$ increases while $E_2$ decreases.

(ii) As $\alpha$ increases ( i.e. if the intercept term $\alpha$ departs from origin in positive direction) relative efficiency of $\bar{y}_a$ with respect to $\bar{y}$ decreases while $E_2$ increases.
(iii) As $\beta$ increases $E_1$ increases for fixed g while $E_2$ is unaffected.
(iv) The maximum gain in efficiency is observed over $\bar{y}$ as well as over $\bar{y}_r$ if A coincide with the value of $\alpha$. Finally, the estimator $\bar{y}_a$ is to be preferred when the intercept term $\alpha$ departs substantially from origin.

Table 1: Relative efficiencies of $\bar{y}_a$ with respect to $\bar{y}$ and $\bar{y}_\Gamma$

| g | β | α = 0.5, n = 10 | | | | | |
|---|---|---|---|---|---|---|---|
| | | $E_1$ | | | $E_2$ | | |
| | | A | | | A | | |
| | | 0.30 | 0.60 | 0.90 | 0.30 | 0.60 | 0.90 |
| 0.0 | 0.5 | 192.86 | 193.23 | 191.40 | 101.34 | 101.54 | 100.57 |
| | 1.0 | 482.16 | 483.16 | 478.09 | 101.34 | 101.54 | 100.57 |
| | 1.5 | 964.32 | 966.17 | 956.98 | 101.34 | 101.54 | 100.57 |
| 0.5 | 0.5 | 132.67 | 132.77 | 132.30 | 100.49 | 100.56 | 100.21 |
| | 1.0 | 237.82 | 237.99 | 237.16 | 100.49 | 100.56 | 100.21 |
| | 1.5 | 413.08 | 413.36 | 411.93 | 100.49 | 100.56 | 100.21 |
| 1.0 | 0.5 | 111.06 | 111.08 | 110.95 | 10.17 | 100.19 | 100.07 |
| | 1.0 | 148.08 | 148.11 | 147.93 | 10.17 | 100.19 | 100.07 |
| | 1.5 | 209.78 | 209.83 | 209.57 | 10.17 | 100.19 | 100.07 |
| 1.5 | 0.5 | 103.99 | 104.00 | 103.96 | 100.06 | 100.07 | 100.03 |
| | 1.0 | 116.64 | 116.65 | 116.60 | 100.06 | 100.07 | 100.03 |
| | 1.5 | 137.71 | 137.72 | 137.66 | 100.06 | 100.07 | 100.03 |
| 2.0 | 0.5 | 102.23 | 102.23 | 102.22 | 100.02 | 100.02 | 100.01 |
| | 1.0 | 106.43 | 106.43 | 106.42 | 100.02 | 100.02 | 100.01 |
| | 1.5 | 113.43 | 113.43 | 113.42 | 100.02 | 100.02 | 100.01 |

| g | β | α = 0.5, n = 20 | | | | | |
|---|---|---|---|---|---|---|---|
| | | $E_1$ | | | $E_2$ | | |
| | | A | | | A | | |
| | | 0.30 | 0.60 | 0.90 | 0.30 | 0.60 | 0.90 |
| 0.0 | 0.5 | 196.58 | 196.96 | 195.11 | 103.33 | 101.52 | 100.56 |
| | 1.0 | 491.46 | 492.39 | 487.77 | 103.33 | 101.52 | 100.56 |
| | 1.5 | 982.92 | 984.39 | 975.53 | 103.33 | 101.52 | 100.56 |
| 0.5 | 0.5 | 134.37 | 134.46 | 134.46 | 100.48 | 100.55 | 100.20 |
| | 1.0 | 240.86 | 241.02 | 240.02 | 100.48 | 100.55 | 100.20 |
| | 1.5 | 418.35 | 418.63 | 417.20 | 100.48 | 100.55 | 100.20 |
| 1.0 | 0.5 | 111.76 | 111.79 | 111.65 | 100.17 | 100.19 | 100.07 |
| | 1.0 | 149.01 | 149.05 | 148.87 | 100.17 | 100.19 | 100.07 |
| | 1.5 | 211.10 | 211.16 | 210.90 | 100.17 | 100.19 | 100.07 |
| 1.5 | 0.5 | 104.00 | 104.00 | 103.96 | 100.06 | 100.07 | 100.02 |
| | 1.0 | 116.64 | 116.65 | 116.60 | 100.06 | 100.07 | 100.02 |
| | 1.5 | 137.71 | 137.73 | 137.67 | 100.06 | 100.07 | 100.02 |
| 2.0 | 0.5 | 101.60 | 101.60 | 101.58 | 100.02 | 100.02 | 100.01 |
| | 1.0 | 105.77 | 105.77 | 105.76 | 100.02 | 100.02 | 100.01 |
| | 1.5 | 112.73 | 112.73 | 112.73 | 100.02 | 100.02 | 100.01 |

Table 2: Relative efficiencies of $\bar{y}_a$ with respect to $\bar{y}$ and $\bar{y}_r$

| | | $\alpha = 1.0$ | | | | | | | |
|---|---|---|---|---|---|---|---|---|---|
| g | $\beta$ | n = 10 | | | | | | | |
| | | $E_1$ A | | | | $E_2$ A | | | |
| | | 0.50 | 1.0 | 1.50 | 1.90 | 0.50 | 1.0 | 1.50 | 1.90 |
| 0.0 | 0.5 | 190.31 | 193.36 | 190.31 | 183.82 | 104.73 | 106.41 | 104.73 | 101.16 |
| | 1.0 | 475.78 | 483.40 | 475.78 | 459.55 | 104.73 | 106.41 | 104.73 | 101.16 |
| | 1.5 | 951.55 | 966.79 | 951.55 | 919.10 | 104.73 | 106.41 | 104.73 | 101.16 |
| 0.5 | 0.5 | 132.03 | 132.80 | 132.03 | 130.34 | 101.73 | 102.32 | 101.73 | 100.43 |
| | 1.0 | 236.67 | 238.05 | 236.67 | 233.65 | 101.73 | 102.32 | 101.73 | 100.43 |
| | 1.5 | 411.07 | 413.46 | 411.07 | 405.82 | 101.73 | 102.32 | 101.73 | 100.43 |
| 1.0 | 0.5 | 110.87 | 111.09 | 110.87 | 110.36 | 100.61 | 100.82 | 100.61 | 100.15 |
| | 1.0 | 147.82 | 148.12 | 147.82 | 147.15 | 100.61 | 100.82 | 100.61 | 100.15 |
| | 1.5 | 209.42 | 209.84 | 209.42 | 208.46 | 100.61 | 100.82 | 100.61 | 100.15 |
| 1.5 | 0.5 | 103.93 | 104.00 | 103.93 | 103.77 | 100.21 | 100.28 | 100.21 | 100.05 |
| | 1.0 | 116.57 | 116.65 | 116.57 | 116.39 | 100.21 | 100.28 | 100.21 | 100.05 |
| | 1.5 | 137.63 | 137.73 | 137.63 | 137.41 | 100.21 | 100.28 | 100.21 | 100.05 |
| 2.0 | 0.5 | 102.21 | 102.23 | 102.21 | 102.15 | 100.67 | 100.09 | 100.07 | 100.01 |
| | 1.0 | 106.41 | 106.43 | 106.41 | 106.3 | 100.67 | 100.09 | 100.07 | 100.01 |
| | 1.5 | 113.41 | 113.43 | 113.41 | 113.35 | 100.67 | 100.09 | 100.07 | 100.01 |

| | | $\alpha = 1.0$ | | | | | | | |
|---|---|---|---|---|---|---|---|---|---|
| g | $\beta$ | n = 20 | | | | | | | |
| | | $E_1$ A | | | | $E_2$ A | | | |
| | | 0.50 | 1.0 | 1.50 | 1.90 | 0.50 | 1.0 | 1.50 | 1.90 |
| 0.0 | 0.5 | 194.01 | 197.08 | 194.01 | 187.47 | 104.67 | 106.33 | 104.67 | 101.14 |
| | 1.0 | 485.03 | 492.70 | 485.03 | 468.68 | 104.67 | 106.33 | 104.67 | 101.14 |
| | 1.5 | 970.06 | 985.40 | 970.06 | 937.36 | 104.67 | 106.33 | 104.67 | 101.14 |
| 0.5 | 0.5 | 133.73 | 134.49 | 133.73 | 132.05 | 101.70 | 102.28 | 101.70 | 100.08 |
| | 1.0 | 239.71 | 241.08 | 239.71 | 236.71 | 101.70 | 102.28 | 101.70 | 100.08 |
| | 1.5 | 416.35 | 418.73 | 416.35 | 411.13 | 101.70 | 102.28 | 101.70 | 100.08 |
| 1.0 | 0.5 | 111.07 | 111.08 | 111.07 | 111.08 | 100.60 | 100.80 | 100.60 | 100.15 |
| | 1.0 | 148.77 | 149.06 | 148.77 | 148.11 | 100.60 | 100.80 | 100.60 | 100.15 |
| | 1.5 | 210.75 | 211.17 | 210.75 | 209.82 | 100.60 | 100.80 | 100.60 | 100.15 |
| 1.5 | 0.5 | 103.94 | 104.01 | 103.94 | 103.78 | 100.20 | 100.27 | 100.20 | 100.05 |
| | 1.0 | 116.57 | 116.65 | 116.57 | 116.40 | 100.20 | 100.27 | 100.20 | 100.05 |
| | 1.5 | 137.64 | 137.73 | 137.64 | 137.42 | 100.20 | 100.27 | 100.20 | 100.05 |
| 2.0 | 0.5 | 101.58 | 101.60 | 101.58 | 101.52 | 100.07 | 100.09 | 100.07 | 100.01 |
| | 1.0 | 105.75 | 105.77 | 105.75 | 105.70 | 100.07 | 100.09 | 100.07 | 100.01 |
| | 1.5 | 112.71 | 112.73 | 112.71 | 112.65 | 100.07 | 100.09 | 100.07 | 100.01 |

Table 3: Relative efficiencies of $\bar{y}_a$ with respect to $\bar{y}$ and $\bar{y}_r$

| | | $\alpha = 1.5$ | | | | | | | | | |
|---|---|---|---|---|---|---|---|---|---|---|---|
| g | $\beta$ | n = 10 | | | | | | | | | |
| | | $E_1$ | | | | | $E_2$ | | | | |
| | | A | | | | | A | | | | |
| | | 0.60 | 1.20 | 1.80 | 2.40 | 2.90 | 0.60 | 1.20 | 1.80 | 2.40 | 2.90 |
| 0.0 | 0.5 | 183.82 | 192.25 | 192.25 | 183.82 | 171.79 | 108.77 | 113.76 | 113.76 | 108.77 | 101.65 |
| | 1.0 | 459.55 | 480.62 | 480.62 | 459.55 | 429.47 | 108.77 | 113.76 | 113.76 | 108.77 | 101.65 |
| | 1.5 | 919.10 | 961.25 | 961.25 | 919.10 | 858.94 | 108.77 | 113.76 | 113.76 | 108.77 | 101.65 |
| 0.5 | 0.5 | 130.34 | 132.52 | 132.52 | 130.34 | 127.01 | 103.29 | 105.01 | 105.01 | 103.29 | 100.64 |
| | 1.0 | 233.64 | 237.55 | 237.55 | 233.65 | 227.67 | 103.29 | 105.01 | 105.01 | 103.29 | 100.64 |
| | 1.5 | 405.82 | 412.60 | 412.60 | 405.82 | 395.44 | 103.29 | 105.01 | 105.01 | 103.29 | 100.64 |
| 1.0 | 0.5 | 110.36 | 111.01 | 111.01 | 110.36 | 109.34 | 101.17 | 101.77 | 101.77 | 101.17 | 100.23 |
| | 1.0 | 147.15 | 148.02 | 148.02 | 147.15 | 147.79 | 101.17 | 101.77 | 101.77 | 101.17 | 100.23 |
| | 1.5 | 208.46 | 209.69 | 209.69 | 208.46 | 206.53 | 101.17 | 101.77 | 101.77 | 101.17 | 100.23 |
| 1.5 | 0.5 | 103.77 | 103.98 | 103.98 | 103.77 | 103.44 | 100.40 | 100.60 | 100.60 | 100.40 | 100.08 |
| | 1.0 | 116.39 | 116.62 | 116.62 | 116.39 | 116.01 | 100.40 | 100.60 | 100.60 | 100.40 | 100.08 |
| | 1.5 | 137.41 | 137.69 | 137.69 | 137.41 | 139.68 | 100.40 | 100.60 | 100.60 | 100.40 | 100.08 |
| 2.0 | 0.5 | 102.15 | 102.22 | 102.22 | 102.15 | 102.04 | 100.13 | 100.20 | 100.20 | 100.13 | 100.03 |
| | 1.0 | 106.35 | 106.42 | 106.42 | 106.35 | 106.24 | 100.13 | 100.20 | 100.20 | 100.13 | 100.03 |
| | 1.5 | 113.35 | 113.42 | 113.42 | 113.35 | 113.23 | 100.13 | 100.20 | 100.20 | 100.13 | 100.03 |

| | | $\alpha = 1.5$ | | | | | | | | | |
|---|---|---|---|---|---|---|---|---|---|---|---|
| G | $\beta$ | n = 20 | | | | | | | | | |
| | | $E_1$ | | | | | $E_2$ | | | | |
| | | A | | | | | | | | | |
| | | 0.60 | 1.20 | 1.80 | 2.40 | 2.90 | 0.60 | 1.20 | 1.80 | 2.40 | 2.90 |
| 0.0 | 0.5 | 187.47 | 196.97 | 195.97 | 187.47 | 175.33 | 108.67 | 113.59 | 113.59 | 108.67 | 101.63 |
| | 1.0 | 468.68 | 489.91 | 489.91 | 468.68 | 438.34 | 108.67 | 113.59 | 113.59 | 108.67 | 101.63 |
| | 1.5 | 937.36 | 979.83 | 979.83 | 937.36 | 876.67 | 108.67 | 113.59 | 113.59 | 108.67 | 101.63 |
| 0.5 | 0.5 | 132.05 | 134.21 | 134.21 | 132.05 | 128.73 | 103.23 | 104.92 | 104.92 | 103.23 | 100.63 |
| | 1.0 | 236.70 | 240.58 | 240.58 | 236.70 | 230.76 | 103.23 | 104.92 | 104.92 | 103.23 | 100.63 |
| | 1.5 | 411.13 | 417.87 | 417.87 | 411.13 | 400.80 | 103.23 | 104.92 | 104.92 | 103.23 | 100.63 |
| 1.0 | 0.5 | 111.08 | 111.72 | 111.72 | 111.08 | 110.08 | 101.14 | 101.72 | 101.72 | 101.14 | 100.23 |
| | 1.0 | 148.11 | 148.96 | 148.96 | 148.11 | 146.77 | 101.14 | 101.72 | 101.72 | 101.14 | 100.23 |
| | 1.5 | 209.82 | 211.02 | 211.02 | 209.82 | 207.92 | 101.14 | 101.72 | 101.72 | 101.14 | 100.23 |
| 1.5 | 0.5 | 103.78 | 103.98 | 103.98 | 103.78 | 103.46 | 100.39 | 100.58 | 100.58 | 100.39 | 100.08 |
| | 1.0 | 116.40 | 116.62 | 116.62 | 116.40 | 116.40 | 100.39 | 100.58 | 100.58 | 100.39 | 100.08 |
| | 1.5 | 137.43 | 137.70 | 137.70 | 137.43 | 137.00 | 100.39 | 100.58 | 100.58 | 100.39 | 100.08 |
| 2.0 | 0.5 | 101.53 | 101.59 | 101.59 | 101.53 | 101.42 | 100.13 | 100.19 | 100.19 | 100.03 | 100.03 |
| | 1.0 | 105.70 | 105.77 | 105.77 | 105.70 | 105.59 | 100.13 | 100.19 | 100.19 | 100.03 | 100.03 |
| | 1.5 | 112.65 | 112.72 | 112.72 | 112.65 | 112.54 | 100.13 | 100.19 | 100.19 | 100.03 | 100.03 |